\documentclass{iucrjournals}

\nolinenumbers

\usepackage{amsmath}
\usepackage{amssymb}
\usepackage{graphicx}   
\usepackage{makecell}
\usepackage{array}
\usepackage{makecell}
\usepackage{endnotes}
\usepackage{caption}
\usepackage{tabularx}
\usepackage{amsmath}
\usepackage{multirow}
\usepackage{booktabs}
\usepackage[justification=centering]{caption}
\usepackage{subcaption}
\title{Tiles from projections of the root and weight lattices of $A_n$}
 
\author[a,*]{Nazife Ozdes Koca}%
\author[a,**]{Mehmet Koca}%
\author[a]{Rehab Nasser Al Reasi}%

\affil[a]{Department of Physics, College of Science, Sultan Qaboos University,  P.O. Box 36, Al-Khoud 123, Muscat, Oman}
\affil[**]{Professor Emeritus}
\affil[*]{e-mail: nazife@squ.edu.om}

\begin{document}
\maketitle
\begin{center}
{\Large \textbf{Abstract}}  
\end{center}
Main purpose of this work is to introduce a general technique of projection of the Voronoi tessellation of the weight lattice $A_n^\ast$ and apply it for the lattice $A_4^\ast$. The projection of the Voronoi tessellation of the weight lattice $A_4^\ast$ produces a totally different tiling scheme than the tiling obtained from the Voronoi cell projection of the lattice $A_4$. The 2D faces of the Voronoi cell of the lattice $A_4^\ast$ are of two types: regular hexagons and squares in 4-dimensions but project into two types of hexagons and two types of rhombuses with edges of two lengths in proportion to golden ratio. The mathematical technique employed is also useful for the projections of the root lattice $A_n$.\newline\noindent\hspace*{0.6cm}
Aperiodic tilings of the plane with $\left(n+1\right)$-fold symmetry can be obtained by projecting the root lattice $A_n$ and the weight lattice $A_n^\ast$ $\left(n=1,\ 2,\ldots\right)$. Projections of the Voronoi cells of these lattices display different tilings of $\left(n+1\right)$-fold symmetry. The orbits of the convex hull of the fundamental weight vectors $\omega_1,\omega_2,\ldots,\omega_n$ under the Coxeter-Weyl group $W\left(a_n\right)$ are the Delone cells and tile the root lattice in an alternating manner as $W{(a_n)(\omega}_1)+W{(a_n)(\omega}_n),W{(a_n)(\omega}_2)+W{(a_n)(\omega}_{n-1}),\ldots,$ where, for example, $W{(a_n)(\omega}_1)+W{(a_n)(\omega}_n)$ \ represents the vector sum of the two orbits. \newline\noindent\hspace*{0.6cm}
A convenient set of linearly dependent and non-orthogonal $\left(n+1\right)$ vectors $k_i$ is introduced 
for the analysis of the polytopes to determine their facets under the group $W\left(a_n\right)$. The simple roots and the fundamental weights are defined as $\alpha_i=k_i-k_{i+1},\left(i=1,2,\ldots,n\right) ,\omega_i=k_1+k_2+\ldots+k_i$, respectively. The point symmetry $W\left(a_n\right) \approx S_{n+1}$ of the lattice permutes the vectors $k_i$. When the vectors $k_i$ are defined in an orthogonal basis, the first two components of $k_i$ determine the Coxeter plane. \newline\noindent\hspace*{0.6cm}
The Delone cell of the weight lattice $A_n^\ast$ is a copy of the fundamental simplex whose vertices are defined as 0, $\omega_1,\omega_2,\ldots,\omega_n$. Projection of the Delone cells of $A_n$ and $A_n^\ast$ on the Coxeter plane displays the same type of tiles and tilings but the Voronoi cell projection of these lattices yields different tiles and tilings. Vertices of the Voronoi cell $V(0)$ of $A_n$ is the union of the orbits of the weight vectors $W(a_n){(\omega}_1)\cup W\left(a_n\right)(\omega_2)\cup\ldots\cup W\left(a_n\right)(\omega_n)$ and the 2D faces are the rhombuses. The Voronoi cell ${V(0)}^\ast$ of $A_n^\ast$ is the permutohedron of order $(n+1)$ and its vertices are the permutations of the vectors ${k}_i$ of the vertex $\frac{1}{n+1}[\left(n+1\right)k_1+nk_2+\ldots+k_{n+1}]$. It has regular hexagons and squares as 2D faces in $n$-dimensions.

\newpage

%%%%%%%%%%%%%%%%%%%%%%%%%%%%%%%%%%%%%%%%%%%%%%%%%%%%%%%%%%

\section{Introduction}
Delone and Voronoi cells of the lattices described by the affine Coxeter-Weyl groups tessellate the $n$-dimensional Euclidean space facet to facet \cite{ConwaySloane1988,ConwaySloane1991}. Projection of higher dimensional lattice is initiated by De Bruijn\cite{deBruijn1981} by projecting the 5D cubic lattice and generalized to the higher dimensional cubic lattices \cite{Whittaker1987}. Study of the facets of the Delone and the Voronoi cells of the root and weight lattices and their projections onto the Coxeter plane are important as they describe the quasicrystallographic structures \cite{Koca2012,Koca2014,Koca2018}. A classification of the affine dihedral subgroups of the affine Coxeter-Weyl groups has been worked out in \cite{Koca2025} including the well-known fact that the affine non-crystallographic groups $W(H_2), W\left(H_3\right)$ and $W(H_4)$ are the respective subgroups of the affine Coxeter-Weyl groups ${W(A}_4), W(D_6)$ and $W(E_8)$. See for further references \cite{PateraTwarock2002} and \cite{Dechant2017}. \newline\noindent\hspace*{0.6cm}
The root lattice $A_n$ and the weight lattice $A_n^\ast$ have interesting quasicrystallographic implications when projected onto the Coxeter plane displaying $\left(n+1\right)$-fold symmetry. Projection of the lattice $A_n$ has been worked out by classifying its Delone cells \cite{Koca2019}. An interesting example is the tessellation displaying 5-fold symmetry obtained by projection of the root lattice $A_4$ tiled by its Delone cells. The Delone cells centralizing the vertices of the Voronoi cell $V(0)$ can be defined by the vector sums of the orbits $W{(a_4)(\omega}_1)+W{(a_4)(\omega}_4),W{(a_4)(\omega}_2)+W{(a_4)(\omega}_3)$. The 2D faces of the Delone cells project as Robinson triangles leading to the Penrose tessellation of the plane with kites and darts. On the other hand, when the Voronoi cell tessellation of $A_4$ is projected onto the Coxeter plane the Penrose tiling is obtained with thick and thin rhombuses. \newline\noindent\hspace*{0.6cm}
Root lattices $A_n$ encode Coxeter-Weyl symmetry and determine invariant projection geometry via the action of the Coxeter element and its associated Coxeter plane, a perspective already visible in the fivefold constructions derived from $A_4$ \cite{Baake1990}. The weight lattices $A_n^\ast$, being dual to $A_n$, frequently govern the geometry of Voronoi cells and hence the structure of acceptance domains and projected prototiles; this explicit connection between the Penrose tiling and projection from $A_4$ or $A_4^\ast$ is emphasized in reference \cite{Heuer2008}. Integrating root and weight duality, the affine Coxeter symmetry, and Voronoi-Delone projections provide a unified framework for extending weight-lattice constructions in aperiodic order.\newline\noindent\hspace*{0.6cm}
The weight lattice $A_4^\ast$ is tiled by the congruent copies of the fundamental simplex determined by the vertices $0,\omega_1,\omega_2,\omega_3,\omega_4$. It is interesting that the 2D faces of the fundamental simplex project as Robinson triangles and leads to the same kind of tiling obtained from the projection of the root lattice $A_4$. However, a more interesting tiling is obtained by projection of the Voronoi tessellation of the $A_4^\ast$. Projection of the Voronoi cell of the weight lattice $A_4^\ast$, the 4-dimensional permutohedron, onto the Coxeter plane yields explicit planar prototiles (thick and thin rhombi together with related hexagonal tiles), as mentioned but not studied in length by \cite{Koca2022}. The main concern of the present work is then a detailed exposition of the Voronoi projection of the lattice $A_4^\ast$. \newline\noindent\hspace*{0.6cm}
In Section 2 we introduce the weight lattice $A_n^\ast$ in its generality by defining the weight vectors  $\omega_i$ in terms of $(n+1)$ linearly dependent non-orthogonal vectors $k_i$ which have convenient representations in an orthogonal basis. The Delone cell is the fundamental simplex and the Voronoi cell is the permutohedron of order $(n+1)$. In Section 3 we discuss the specific properties of the permutohedron of order 5, the Voronoi cell of $A_4^\ast$. In Section 4 we project the 2D faces of the permutohedron of order 5 which are regular hexagons and squares in 4D and prove that they project onto four different tiles; two rhombuses and two hexagons with edges proportional to $1$ and the golden ratio $\tau$. A tessellation of the plane with these four tiles has been introduced. A discussion of our result and its generalization to the Voronoi projection of the weight lattice $A_n^\ast$ for $n\geq5$ has been presented.

%%%%%%%%%%%%%%%%%%%%%%%%%%%%%%%%%%%%%%%%%%%%%%%%%%%%%%

\section{The weight lattice \texorpdfstring{$A_n^\ast$}{}}
The Coxeter group ${W(a}_n)$ is defined by its reflection generators as,
\begin{equation}
W\left(a_n\right):=\left\langle r_1,r_2,\ldots,r_n\middle|\left(r_ir_j\right)^{m_{ij}}=1\right\rangle,
\end{equation}
where $r_i$ is the reflection generator associated with the simple root $\alpha_i=k_i-k_{i+1}$, which acts on an arbitrary vector $\lambda$ as $r_i\lambda=\lambda-\frac{\left.2(\lambda,\alpha_i\right)}{\left(\alpha_i, \alpha_i\right)}\alpha_i$ leading to $r_i$: $k_i\leftrightarrow k_{i+1}$. The non-orthogonal and linearly dependent vectors $k_i,(i=1,\ldots,n+1)$ can be defined in terms of linearly independent and orthonormal set of vectors $l_i,\ (i=1,\ldots,n+1)$ as $k_i=-l_i+\frac{l_1+\ldots+l_{n+1}}{n+1}$. Here $m_{ii}=1 (i=1,\ldots,n), m_{i,i+1}=m_{i+1,i}=3 (i=1,\ldots,n-1),$ and $m_{ij}=m_{ji}=2\ (1\le i\le j-1\le n-1)$. 
The Coxeter-Dynkin diagram given in Fig. 1 explains all properties of the Coxeter group \cite{CoxeterMoser1965,Coxeter1973,Humphreys1992}   
\begin{figure}[ht]
\begin{center}
\includegraphics[width=0.3\textwidth]{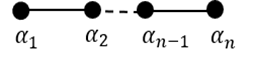} 
\end{center}
\caption{Coxeter-Dynkin diagram $a_n$.}
\end{figure}
\newline
The root system and certain polytopes of the group ${W(a}_n)\approx S_{n+1}$ has a symmetry called Dynkin diagram symmetry $\gamma$ which acts on the simple roots and weights as 
\begin{equation}
\gamma: \alpha_1\longleftrightarrow\alpha_n,\ \alpha_2\longleftrightarrow\alpha_{n-1},\ldots;   \ \ \ \gamma: \omega_1\longleftrightarrow\omega_n,\ \omega_2\longleftrightarrow\omega_{n-1},\ldots.           
\end{equation}
Therefore, the root system and certain polytopes symmetric under the Dynkin diagram symmetry has a larger symmetry ${W(a}_n):C_2$ of order $2\times\left(n+1\right)!$. \newline
The weights $\omega_i$ are defined by the inner product $\left(\alpha_i,\omega_j\right)=\delta_{ij}$ satisfying the relations $\omega_i=\sum_{j}{({C^{-1})}_{ij}\alpha_j,} \  \alpha_i=\sum_{j} {C_{ij}\omega_j,} \ (\omega_i,\omega_j)=({C^{-1})}_{ij}$ where $C_{ij}$ is the matrix element of the Cartan matrix. The weight lattice $A_n^\ast$ is defined as a linear combination of the weight vectors with integer coefficients
\begin{equation}
\Lambda^\ast: q=\sum_{i=1}^{n}n_i\omega_{i\ \ \ \ },\ \ n_i\in\mathbb{Z}, \ \ \ \omega_i=k_1+k_2+\ldots+k_i,      
\end{equation}
where $k_j$ is given in an orthogonal basis as
\begin{equation}
\begin{aligned}
k_j &= \sqrt{\frac{2}{n+1}}\big(
\cos\left(\frac{2\pi j}{n+1}\right),
\sin\left(\frac{2\pi j}{n+1}\right), \\
&\qquad
\cos\left(\frac{4\pi j}{n+1}\right),
\sin\left(\frac{4\pi j}{n+1}\right),
\ldots
\big), \ \ (j=1,2,\ldots,n+1) \ (\text{for even $n$}).
\end{aligned}
\end{equation}
For odd $n$ add another component $\frac{{(-1)}^j}{\sqrt2}$. The simple roots of the Lie algebra of $a_n$ are given by $\alpha_i=k_i-k_{i+1}$ and the root lattice $A_{n}$ is defined as the set of linear combinations of the simple root vectors $\Lambda: p=\sum_{i=1}^{n}{m_i\alpha_i},{\ m}_i\in\mathbb{Z}$. Note that the root and weight lattices are the reciprocal lattices for the inner product of the vectors $\left(p,q\right)=$ integer. The vectors $k_j$ are special in the sense that they represent the vertices of $n$-simplex. For example, for $n=3$ they represent the vertices of a tetrahedron and for $n=4$ they are the vertices of a 4-simplex. The vectors $k_i$ satisfy the relations
\begin{equation}
k_1+k_2+\ldots+k_{n+1}=0,\ \ \left(k_i,k_i\right)=\frac{n}{n+1}, \ \ \left(k_i,k_j\right)=-\frac{1}{n+1}\ ,\ \ (i\neq j=1,\ 2,\ldots,n+1).         
\end{equation}
The equations (4-5) imply that we work in $n$-dimensional Euclidean space.  \newline\noindent\hspace*{0.6cm}
The vectors $k_j$ in (4) are defined in such a way that the first two components generate the Coxeter plane. Delone cells of the root lattice $A_n$ have been already discussed \cite{Koca2019} and the projections of their 2D faces (equilateral triangles in higher dimensions) are classified as triangles; Robinson triangles for $n=4$ and Danzer triangles \cite{NischkeDanzer1996} for $n=6$ and many more triangular tiles for other higher values of n have been worked out. Voronoi cell projections of the lattice $A_n$ have also been classified \cite{Koca2019} and shown that they result in various rhombuses; thick and thin Penrose rhombuses for $n=4$. \newline\noindent\hspace*{0.6cm}
The projection of the weight lattice $A_n^\ast$ has not been studied in its generality via its Delone cells and its Voronoi cell. Preliminary study implies that the quasi-crystallographic structures arising from the Voronoi cell projection offers different perspectives as we will discuss in Section 3 for the lattice $A_4^\ast$. \newline\noindent\hspace*{0.6cm}
Now, first let us discuss the Delone cell of the lattice $A_{n}^\ast$ and its projection. The fundamental simplex with $\left(n+1\right)$ vertices $0,\omega_1,\omega_2,\ldots,\omega_n$ \cite{ConwaySloane1991} is the Delone cell of the lattice $A_{n}^\ast$ which tiles the lattice. For example, the vertices $0,\omega_1,\omega_2$ form an equilateral triangle and constitute the Delone cell of the self-dual lattice $A_2^\ast\approx A_2$. For $n=3$, the Delone cell $0,\omega_1,\omega_2,\omega_3$ is an irregular tetrahedron covered by two different isosceles triangular faces. One of the vertices of the Voronoi cell ${V(0)}^\ast$ of the lattice $A_{n}^\ast$ is obtained by averaging the vertices of the fundamental simplex defined as
\begin{equation}
q=\frac{1}{n+1}[0+\omega_1+\omega_2+\ldots+\omega_n], 
\end{equation}
which can also be written as
\begin{equation}
q=\frac{1}{n+1}[\left(n+1\right)k_1+nk_2+\ldots+k_{n+1}], 
\end{equation}
\noindent\hspace*{0.6cm}
The Voronoi cell centered around the origin is a permutohedron with $\left(n+1\right)!$ vertices and the Coxeter-Weyl group $W(a_n) \approx S_{n+1}$ permutes the $\left(n+1 \right)$ vectors $k_j$ in (7). Permutohedron is invariant under the Dynkin diagram symmetry implying that it has the symmetry of group of order $2\times\left(n+1\right)!$. Permutation of the vectors $k_j$ induces the permutations of the integer coefficients $\left(n+1\right), n,\ldots,1$. Therefore, the vertex in (7) other than the term in the denominator can be symbolically represented as a sequence of the integers $[\left(n+1\right)n(n-1)\ldots2\ 1]$. Then the permutohedron is the mathematical abstraction of the vertex (7) excluding the term in the denominator \cite{Ziegler1995}. \newline\noindent\hspace*{0.6cm}
For $n=2$ the Voronoi cell is a regular hexagon and for $n=3$ the Voronoi cell of the lattice $A_3^\ast$ is a truncated octahedron consisting of hexagonal and square faces. For arbitrary value $n$ the Voronoi cell is a permutohedron of order $(n+1)$ tiled by the $\left(n+1\right)!$ fundamental simplices. The permutohedron has various facets which will be discussed for $n=4$ in details in Section 3. \newline\noindent\hspace*{0.6cm}
Projection of the vertices of the Delone cell of the lattice $A_n^\ast$ onto the Coxeter plane form a regular polygon with $n+1$ vertices. Projections of the 2D faces are then the triangles formed by any three vertices. Classification of these triangles for arbitrary $n$ follows the same pattern obtained by the Delone cell projection of the lattice $A_n$ which has been worked out in the reference \cite{Koca2019}. \newline\noindent\hspace*{0.6cm}
The Voronoi cell projection of the weight lattice $A_n^\ast$ is more intriguing because it involves huge numbers of hexagons and squares as $n$ increases. Their projections and the emerging tessellation should be worked out with care. In the next section we will discuss the Voronoi cell projection of the lattice $A_4^\ast$. We are motivated from the observation that the Voronoi cell of $A_4^\ast$ is a permutohedron of order 5 and its projection results in four different tiles as noted by \cite{Koca2022} but not discussed in details. In what follows we present a detailed account of the Voronoi tessellation of the lattice $A_4^\ast$ and its projection onto the Coxeter plane.

%%%%%%%%%%%%%%%%%%%%%%%%%%%%%%%%%%%%%%%%%%%%%%%%%%%%%%%%

\section{Voronoi cell of the weight lattice \texorpdfstring{$A_4^\ast$}{}}
In this section, we study the structure of the facets of the Voronoi cell of $A_4^\ast$. Center of the Delone cell with the vertices $0,\omega_1,\omega_2,\omega_3,\omega_4$ is the vector
\begin{equation}
q=\frac{1}{5}[5k_1+4k_2+3k_3+2k_4+k_5],
\end{equation}
where 
\begin{equation}
\left[5k_1+4k_2+3k_3+2k_4+k_5\right]=[\omega_1+\omega_2+\omega_3+\omega_4]
\end{equation}
can be denoted as the sequence of integers 54321 as will be explained below. \newline\noindent\hspace*{0.6cm}
Permutation of the vectors $k_i$ in (9) under the group $W(a_4)\approx S_5$ forms a 4D polytope with 120 vertices \cite{Koca2012}. This polytope is also known as the permutohedron of order 5 whose symmetry is the semi-direct product of the Dynkin diagram symmetry with the Coxeter group $S_5:C_2\approx W\left(a_4\right):C_2= \langle r_1, r_2, r_3, r_4 \rangle:C_2$.  \newline\noindent\hspace*{0.6cm}
The group $W\left(a_4\right)$ acting on the vector (9) induces a permutation on the sequence of integers 54321. The polytope has the following facets:
\begin{equation}
\begin{aligned}
N_0&=120\ \text{vertices,} \\
N_1&=240\ \text{edges}, \\
N_2&=150\ (60\ \text{hexagons} +90\ \text{squares}), \\
N_3&=30\ (10\ \text{truncated octahedra} +20\ \text{hexagonal prisms}).
\end{aligned}
\end{equation}
\noindent\hspace*{0.6cm}
These numbers satisfy the Euler characteristic $N_0-N_1+N_2-N_3=0$. Note that it is important to classify the 2D faces of the Voronoi cells because their projections will determine the structures of the tiles. Classification of the 2D faces can be well understood by studying the 3D faces of the Voronoi cell. \newline\noindent\hspace*{0.6cm}
To determine the facets of the permutohedron, one uses the coset decomposition of the group $W\left(a_4\right)$ under the subgroup which determines the symmetry of the related facet. The symmetries of some particular 3D facets can be represented by the groups $\langle r_1, r_2, r_3 \rangle$ and  $\langle r_2, r_3, r_4 \rangle$ of orders 24 each describing the symmetry of a truncated octahedron so that the number of truncated octahedral facets is $10=5+5$. These two groups are conjugate to each other under the Dynkin diagram symmetry $\gamma$. Another 3D facet is the hexagonal prism either described by the symmetry $\langle r_1, r_2, r_4 \rangle$ or $\langle r_1, r_3, r_4 \rangle$, both of order 12, which are also $\gamma$ conjugate of each other and hence the number of hexagonal prisms can be written as $20=10+10$.\newline \newline
\textbf{Determination of the hexagonal faces as well as the centers of the truncated octahedra} \newline \newline
To obtain the vertices of one of the truncated octahedron, first apply the hexagonal symmetry $\langle r_1, r_2 \rangle$ on the vertex [54321] horizontally and then apply a  group element $a=r_1r_2r_3$ on each vertex of hexagons, which permutes the first 4 integers in the cyclic order but fixes the integer 1.  Then, one obtains the list of vertices in Table 1.
\begin{table}[h]
\centering
\caption{List of vertices of a truncated octahedron with its 4 hexagonal faces.}
\begin{tabular}{c c c c c c c c}
\toprule
 & \multicolumn{6}{c}{Apply the group $\langle r_1, r_2 \rangle$} & \multirow{2}{*}{Center of each hexagon} \\
\cmidrule(lr){2-7}
Apply $a=r_1 r_2 r_3$ 
& $r_1$ & $r_1 r_2 r_1$ & $r_2$ & $r_1$ & $r_1 r_2 r_1$ &  &  \\
$\downarrow$ 
& $\rightarrow$ & $\rightarrow$ & $\rightarrow$ & $\rightarrow$ & $\rightarrow$ &  &  \\
\midrule

$1^{\text{st}}$ hexagon 
& 54321 & 45321 & 35421 & 34521 & 43521 & 53421 
& $-2k_4 - 3k_5$ \\

$2^{\text{nd}}$ hexagon 
& 25431 & 24531 & 23541 & 23451 & 24351 & 25341 
& $-2k_1 - 3k_5$ \\

$3^{\text{rd}}$ hexagon 
& 32541 & 32451 & 42351 & 52341 & 52431 & 42531 
& $-2k_2 - 3k_5$ \\

$4^{\text{th}}$ hexagon 
& 43251 & 53241 & 54231 & 45231 & 35241 & 34251 
& $-2k_3 - 3k_5$ \\

\bottomrule
\end{tabular}
\end{table}
\newline
Table 1 shows the centers of 4 hexagons of the truncated octahedral facet of the permutohedron whose center is  $-\frac{5}{2}k_5$. Considering the actual Voronoi cell one should divide this by $5$, then the center of this facet is $-\frac{1}{2}k_5$. Applying a group element, say, $b=r_1r_2r_3r_4$ of order $5$ one creates altogether 5 truncated octahedra with centers $-\frac{1}{2}k_j, \left(j=1,2,3,4,5\right)$. Using the Dynkin diagram symmetry $\gamma:k_1\longleftrightarrow-k_5,k_2\longleftrightarrow-k_4,k_3\rightarrow-k_3$, the centers of $10$ truncated octahedral facets can be written as 
\begin{equation}
\begin{aligned}
\pm\frac{1}{2}k_j, \  \left(j=1,2,3,4,5\right).    
\end{aligned}
\end{equation}
The centers of the 20 hexagonal faces of 5 truncated octahedra obtained by the coset decomposition of the group $W\left(a_4\right)$ under the hexagonal symmetry $\langle r_1, r_2 \rangle$ are the copies of the centers of hexagons in Table 1. Permuting the centers of 4 hexagons in Table 1 by the group element b, one finds the centers of 20 hexagonal faces as $-\frac{1}{5}{(2k}_i+3k_j), \left(i\neq j=1,2,3,4,5\right)$. They can also be written as the orbit $W\left(a_4\right)\frac{1}{5}\left(2\omega_3+\omega_4\right)$. \newline
The centers of the next 20 hexagons can be obtained by applying the Dynkin diagram symmetry $\gamma : W\left(a_4\right)\frac{1}{5}\left(2\omega_3+\omega_4\right)= W\left(a_4\right)\frac{1}{5}\left(2\omega_2+\omega_1\right).$  \newline
The centers of the remaining 20 hexagons can be determined by the coset decomposition of the group $W\left(a_4\right)$ under the hexagonal symmetry $\langle r_2, r_3 \rangle$ which is invariant under the Dynkin diagram symmetry $\gamma$. The hexagon generated by the group $\langle r_2, r_3 \rangle$ acting on the vertex $\left[54321\right]$ has the center $\frac{2}{5}{(k}_1-k_5)$. Application of the generator $a=r_1r_2r_3$ on the vertices of this hexagon reorders the vertices in Table 1 and leads to the centers of hexagons as $\frac{2}{5}{(k}_i-k_5), \left(i=1,2,3,4\right)$. Applying the group element $b=r_1r_2r_3r_4$ on the vector $\frac{2}{5}{(k}_i-k_5)$ we obtain the third set of 20 hexagons with centers represented by the vectors $\frac{2}{5}{(k}_i-k_j),(i\neq j=1,2,3,4,5)$. These are the vectors proportional to the set of roots. This is a self-dual orbit $W\left(a_4\right)\frac{2}{5}\left(\omega_1+\omega_4\right)$ under the Dynkin diagram symmetry. As we have pointed out in (10) the number of hexagons can be written as $60=20+20+20$, two sets of which are conjugates of each other under the symmetry $\gamma$ and the third one is the self-conjugate. The centers of the hexagons are represented by the following orbits:
\begin{equation}
\begin{aligned}
W\left(a_4\right)\frac{1}{5}\left(2\omega_3+\omega_4\right), W\left(a_4\right)\frac{1}{5}\left(2\omega_2+\omega_1\right),W\left(a_4\right)\frac{2}{5}\left(\omega_1+\omega_4\right).
\end{aligned}
\end{equation}
\newline
\textbf{Centers of the hexagonal prisms and the squares} \newline \newline
A hexagonal prism is determined by applying an orthogonal generator, say $r_4$, on the vertices of a hexagon determined by the orbit $\langle r_1, r_2 \rangle \left[54321\right]$ which is shown in Table 2.
\begin{table}[h]
\centering
\caption{Vertices of a hexagonal prism.}
\begin{tabular}{c c c c c c c c c c c c}
\toprule
 & \multicolumn{11}{c}{Apply the group $\langle r_1, r_2 \rangle$} \\
\cmidrule(lr){2-12}
 & & $r_1$ &  & $r_1 r_2 r_1$ &  & $r_2$ &  & $r_1$ &  & $r_1 r_2 r_1$ &  \\
& & $\rightarrow$ &  & $\rightarrow$ &  & $\rightarrow$ &  & $\rightarrow$ &  & $\rightarrow$ &  \\
\midrule
Apply $r_4 \downarrow$
& 54321 &  & 45321 &   & 35421 &   & 34521 &  & 43521 &   & 53421 \\
& 54312 &  & 45312 &  & 35412 &  & 34512 &  & 43512 & & 53412 \\
\bottomrule
\end{tabular}
\end{table}
\newline
The center of this prism is determined by averaging its 12 vertices as $\frac{1}{2}\omega_3=-\frac{1}{2}{(k}_4+k_5)$. The orbit $W\left(a_4\right)\frac{1}{2}\omega_3$ consists of 10 vectors $-\frac{1}{2}{(k}_i+k_j),\ (i\neq j=1,2,3,4,5)$. The conjugate of this set of vectors, under the Dynkin diagram symmetry, is the orbit $W\left(a_4\right)\frac{1}{2}\omega_2$ and represents the  centers of the other 10 hexagonal prisms. Altogether the centers of the hexagonal prisms can be written as  
\begin{equation}
\begin{aligned}
\pm\frac{1}{2}{(k}_i+k_j),\ (i\neq j=1,2,3,4,5).   
\end{aligned}
\end{equation}
Total number of square faces of the permutohedron of order 5 is $90=30+30+30$ obtained as the coset decompositions of $W\left(a_4\right)$ under the groups $\langle r_1, r_4 \rangle, \langle r_1, r_3 \rangle$ and $\langle r_2, r_4 \rangle$. Here the group $\langle r_1, r_4 \rangle$ is self-conjugate and the next two groups are conjugates of each other under the Dynkin diagram symmetry $\gamma$. Centers of the squares can be determined as the orbits
\begin{equation}
W\left(a_4\right)\frac{3}{10}\left(\omega_2+\omega_3\right),  W\left(a_4\right)\frac{1}{10}\left(4\omega_2+3\omega_4\right),  W\left(a_4\right)\frac{1}{10}\left(3\omega_1+4\omega_3\right).
\end{equation}

%%%%%%%%%%%%%%%%%%%%%%%%%%%%%%%%%%%%%%%%%%%%%%%%%%%%%%

\section{Projections of 2D faces of the Voronoi cell of the weight lattice \texorpdfstring{$A_4^\ast$}{} and tilings}
Vertices of a 4-simplex $(k_1,\ldots,\ k_5)$ can be written in an orthogonal basis as 
\begin{equation*}
k_j=\sqrt{\frac{2}{5}}\left(\cos{\frac{2\pi j}{5}},\sin{\frac{2\pi j}{5}},\cos{\frac{4\pi j}{5}},\sin{\frac{4\pi j}{5}}\right), \left(j=1,2,3,4,5\right),
\end{equation*}
or explicitly,
\begin{equation}
\begin{aligned}
k_1&=\frac{1}{\sqrt{10}}\left(-\sigma,\sqrt{2+\tau},-\tau,\sqrt{2+\sigma}\right), \ \ \ k_2=\frac{1}{\sqrt{10}}\left(-\tau,\sqrt{2+\sigma},-\sigma,-\sqrt{2+\tau}\right), \\
k_3&=\frac{1}{\sqrt{10}}\left(-\tau,-\sqrt{2+\sigma},-\sigma,\sqrt{2+\tau}\right), \ \  k_4=\frac{1}{\sqrt{10}}\left(-\sigma,-\sqrt{2+\tau},-\tau,-\sqrt{2+\sigma}\right), \\
k_5&=\frac{1}{\sqrt{10}}\left(2,\ 0,\ 2,\ 0\right), \  \text{with} \ \ \tau=\frac{1+\sqrt5}{2} , \ \    \sigma=\frac{1-\sqrt5}{2}.            
\end{aligned}
\end{equation}
They satisfy the inner products 
\begin{equation}
\begin{aligned}
\left(k_i,k_i\right)=\frac{4}{5}, \ \left(k_i,k_j\right)=-\frac{1}{5}\ ,\ \ \left(i\neq j=1,2,3,4,5\right),        
\end{aligned}
\end{equation}
with the magnitude $\parallel k_i-k_j\parallel=\sqrt2$. \newline
Vertices of the 4-simplex lie on a 3-sphere with radius $\frac{2}{\sqrt{5}}$. Projection of the 4-simplex onto the Coxeter plane form a circle of radius $\sqrt{\frac{2}{5}}$ determined by the first two components of the vectors $k_j$. The magnitudes of the projected vectors $ (k_i-k_j)$ are $\sqrt{\frac{2}{2+\tau}}$ if $k_i$ and  $k_j$ are adjacent to each other on the circle, otherwise, they are equal to $\tau\sqrt{\frac{2}{2+\tau}}$. If the overall factor $\sqrt{\frac{2}{2+\tau}}$ is ignored they are represented as $1$ and $\tau$, respectively. \newpage
\textbf{Projections of hexagons} \newline
The 60 hexagons described above can be classified in two classes represented by the edges determined by the following sets of vectors
\begin{equation}
\begin{aligned}
A:  k_2-k_1,\ {\ k}_3-k_1,\ {\ k}_3-k_2,\ {\ k}_1-k_2,\ {\ k}_1-k_3,\ \ k_2-k_3;  \\  
B:  k_2-k_1,\ {\ k}_4-k_1,\ {\ k}_4-k_2,\ \ k_1-k_2,\ {\ k}_1-k_4,\ {\ k}_2-k_4.    
\end{aligned}
\end{equation}
The first set of edges of the hexagon is obtained by applying the group $\langle r_1, r_2 \rangle$ on the vector $\left[54321\right]$. When we apply the generators of $W\left(a_4\right)$ say, $r_4$, on the set of vertices in the set A the set does not change because $r_4$ commutes with $r_1$ and $r_2$. Applying $r_1$\ and $r_2$ does not change the set of vertices in A but reorders the edges. The generator $r_3$ exchanges the sets A and  B. This proves that 60 hexagons are classified as the sets A and B under the group $ W\left(a_4\right)$. We also note that the Dynkin diagram symmetry does not change the sets. Therefore, it is sufficient to work with the projections of two hexagons represented by the sets A and B. \newline\noindent\hspace*{0.6cm}
When the sets of hexagons are projected onto the Coxeter plane they will be represented by the thin and thick hexagons whose edge lengths are proportional to $\left(1,\tau,1,1,\tau,1\right)$ and $\left(1,\tau,\tau,1,\tau,\tau\right)$, respectively, as shown in Fig. 2.
\begin{figure}[ht]
\begin{center}
\includegraphics[width=0.3\textwidth]{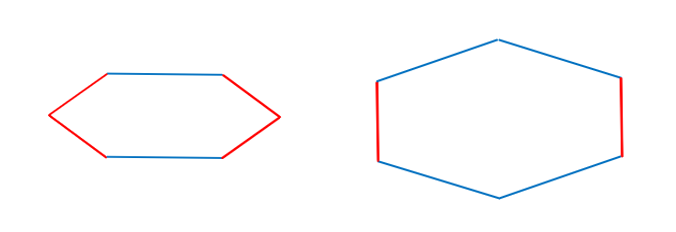}
\caption{Thin and thick hexagons projected from the sets A and B.}
\end{center}
\end{figure} 
\newline
\textbf{Projections of squares} \newline
Square faces in 4D are classified as three sets of squares with edges similar to:
\begin{equation}
\begin{aligned}
C:&  \pm\left(k_1-k_2\right),\pm\left({\ k}_4-k_5\right); \\
D: & \pm\left(k_1-k_3\right),\pm\left({\ k}_{\ 2}-k_4\right); \\
E: &\pm\left(k_1-k_2\right),\pm\left({\ k}_3-k_5\right).
\end{aligned}
\end{equation}
The edges of the square in the set C are obtained from the group action $\langle r_1, r_4 \rangle$ on the vector $\left[54321\right]$ and it is self-conjugate under the Dynkin diagram symmetry. The next two sets of squares are represented by the sets D and E. \newline\noindent\hspace*{0.6cm}
The square represented by the set C projects onto a thin rhombus with interior angles 36° and 144° with edges proportional to 1. The square represented by the set D projects onto a thick rhombus with interior angles 72° and 108° with edges proportional to the golden ratio $\tau$. The square in the set E projects onto a line segment. The thin (red) and thick (blue) rhombuses with edges $1$ and $\tau$ are depicted in Fig. 3.
\begin{figure}[ht]
\begin{center}
\includegraphics[width=0.25\textwidth]{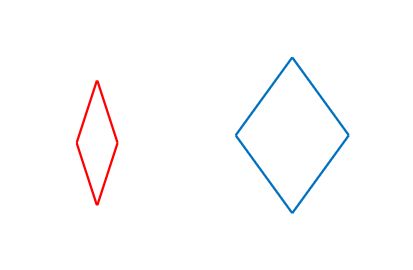}
\caption{Thin (red) and thick (blue) rhombuses with respective edges $1$ and $\tau$.}
\end{center}
\end{figure} 
\newline
We have proved that, by means of the group theoretical analysis, the 2D faces of the Voronoi cell of the weight lattice $A_4^\ast$ project into 4 different tiles given in Fig. 2 and Fig. 3. The projected Voronoi cell, possessing only one reflection symmetry, is tiled by these four types of tiles, five from each type, as shown in Fig. 4.
\begin{figure}[ht]
\begin{center}
\includegraphics[width=0.25\textwidth]{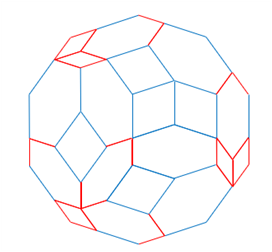}
\caption{Tiling of the projected Voronoi cell of $A_4^\ast$ by the four tiles.}
\end{center}
\end{figure}
\newline
Two patches of centrally symmetric tilings with the four tiles are illustrated in Fig. 5 with uncolored (a) and colored tiles (b).  
\begin{figure}[ht]
\centering
\begin{subfigure}{0.25\textwidth}
\centering
\includegraphics[width=\linewidth]{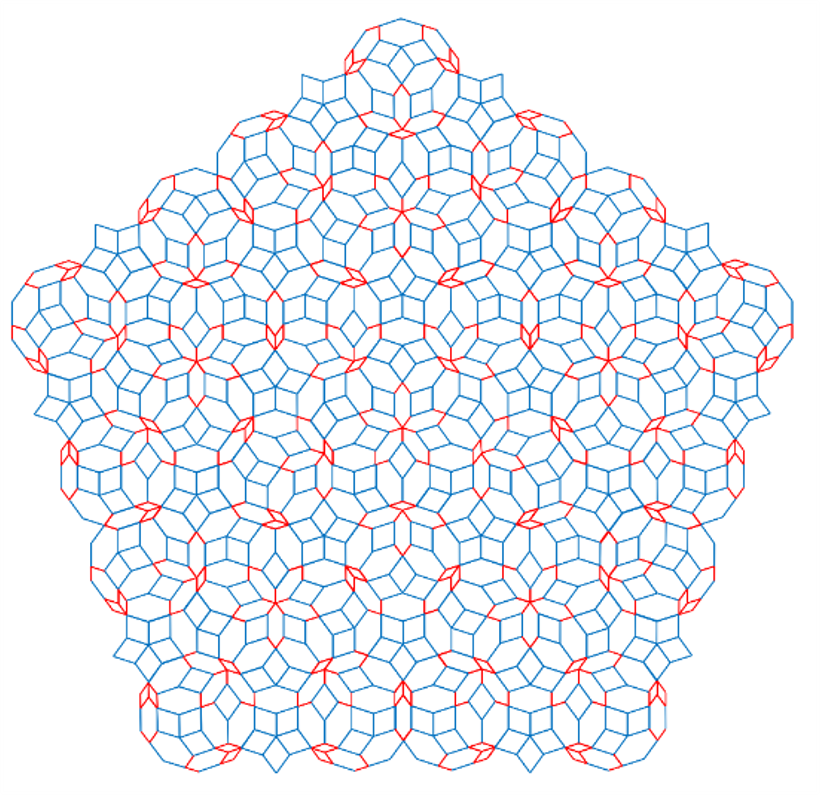}
\caption{}
\end{subfigure}
\hspace{0.02\textwidth}
\begin{subfigure}{0.25\textwidth}
\centering
\includegraphics[width=\linewidth]{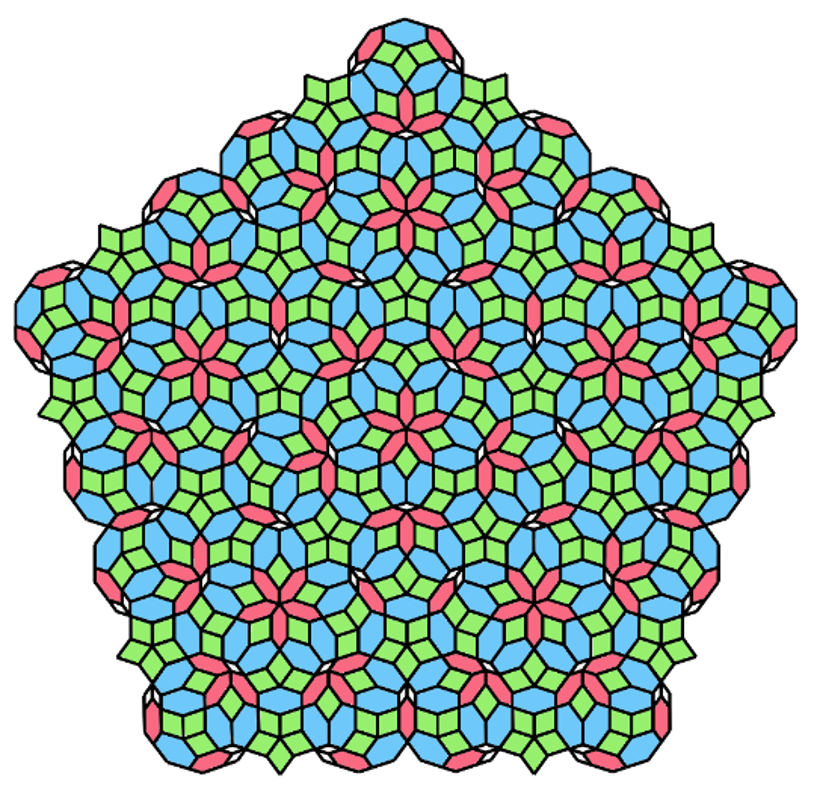}
\caption{}
\end{subfigure}
\caption{Centrally symmetric uncolored (a) and colored (b) patches obtained by projection of the Voronoi cell tessellation of $A_4^\ast$.}
\end{figure}
\newline

%%%%%%%%%%%%%%%%%%%%%%%%%%%%%%%%%%%%%%%%%%%%%%%%%%%%%%%

{\Large \textbf{Conclusion}}  \newline\newline
A general exposition on the projections of the lattices $A_n$ and $A_n^\ast$ via their Delone cells and their Voronoi cells has been introduced. Delone cells tiling the lattice $A_n$ and the lattice $A_n^\ast$ project onto exactly similar triangles which are the Robinson triangles for $n=4$ and the Danzer triangles for $n=6$ and so on. In this work, we emphasized on the structure of the Voronoi polytope of $A_4^\ast$ and its projection. The projection of the Voronoi cell of $A_4$ leads to the Penrose thick and thin rhombuses, on the other hand,  the projection of the Voronoi cell of $A_4^\ast$ consists of four types of tiles; two of them are irregular hexagons with edges proportional to $(\tau,1,1,\tau,1,1)$ and $(\tau,\tau,1,\tau,\tau\,1)$ and a thin and a thick rhombuses with edge lengths $(1,1,1,1)$ and $(\tau,\tau\,\tau,\tau)$, respectively. A centrally symmetric tessellation of the Coxeter plane has been obtained. The study of the permutohedra, the Voronoi cells of the weight lattice $A_n^\ast$, for $n\geq5$ is under investigation and their projection may lead to interesting results. The projection of the Voronoi cell for $n=7$ and $n=11$ can be useful for the 8-fold and 12-fold symmetric quasicrystallographic tilings.
\newline \newline
{\Large \textbf{Acknowledgment}}  \newline\newline
We would like to thank Professor Ramazan Koc for helpful discussions and the referees for constructive and insightful comments.

\bibliography{iucr} 

\end{document}